\documentclass{amsart}
\usepackage{verbatim}
\usepackage{amsmath}
\usepackage{amsthm}
\usepackage{amssymb}
\usepackage{setspace}
\usepackage{xargs}[2008/03/08]
\usepackage{microtype}

\makeatletter
\theoremstyle{plain}
\newtheorem{thm}{\protect\theoremname}[section]
\theoremstyle{definition}
\newtheorem{defn}[thm]{\protect\definitionname}
\theoremstyle{plain}
\newtheorem{lem}[thm]{\protect\lemmaname}
\theoremstyle{remark}
\newtheorem{rem}[thm]{\protect\remarkname}

\makeatother

\providecommand{\definitionname}{Definition}
\providecommand{\lemmaname}{Lemma}
\providecommand{\remarkname}{Remark}
\providecommand{\theoremname}{Theorem}
\let\cal=\mathcal

\begin{document}

\title{Are zero-symmetric simple nearrings with identity equiprime?}

\author{Wen-Fong Ke}

\address{Department of Mathematics, National Cheng Kung University,
Tainan 70101, Taiwan}
\email{wfke@mail.ncku.edu.tw}

\author{Johannes H. Meyer}

\address{Department of Mathematics and Applied Mathematics, University of the Free State, PO Box 339, Bloemfontein, 9300, South Africa}
\email{meyerjh@ufs.ac.za}

\begin{abstract}
We show that there exist zero-symmetric simple nearrings with identity
which are not equiprime, solving a longstanding open problem.
\end{abstract}



\maketitle

\section{Introduction}

\global\long\def\End{\textup{End}}%
\newcommandx\gen[1][usedefault, addprefix=\global, 1=]{\langle#1\rangle}%
\global\long\def\Im{\textup{Im}}%
\global\long\def\Sym{\textup{Sym}}%
\global\long\def\E{\Sigma}%
\global\long\def\N{\mathbb{N}}%
The question stated in the title appears in the 1990 paper \cite{BoothGV}
by Booth, Groenewald, and Veldsman. We answer this question negatively
by showing that certain zero-symmetric simple nearrings $(N,+,\cdot)$
with identity, obtained in another work \cite{KeMPW}, are not equiprime.
The additive groups of these nearrings are infinite simple groups
constructed using the Higman-Neumann-Neumann (HNN) extensions (cf.~\cite{HigmanNN49}).

Here, a \textit{nearring} (actually a \textit{right nearring}) is
a triple $(N,+,\cdot)$ where $(N,+)$ is a group (with neutral element~$0$),
$(N,\cdot)$ is a semigroup, and the right distributive law holds,
i.e., for all $a,b,c\in N$, $(a+b)\cdot c=a\cdot c+b\cdot c$. The
reader is referred to \cite{Pilz83} for general definitions and results
with respect to the theory of nearrings.

Now, a nearring $N$ is \textit{zero-symmetric} if $0a=a0=0$ for
all $a\in N$. It has an identity if there exists an element $1\in N$
such that $1a=a1=a$ for all $a$ in $N$. Additionally, $N$ is \textit{simple}
if its only ideals are the zero ideal and the entire nearring $N$
itself. Analogous to rings, an ideal of $N$ corresponds to the kernel
of a nearring homomorphism on $N$, and thus, $N$ is simple if all
nontrivial homomorphisms of $N$ are monomorphisms.

Unlike the ring case, there are several variants of primeness defined
for nearrings. For instance, in \cite{Veldsman}, $k$-prime nearrings,
$k\in\{0,1,2,3\}$, are defined. However, as Kaarli and Kriis demonstrated
in \cite{KaaKri}, none of these definitions leads to a Kurosh-Amitsur
prime radical class in the variety of zero-symmetric nearrings. This
result prompted the paper by Booth, Groenewald, and Veldsman \cite{BoothGV},
where the notion of an equiprime nearring was introduced.
\begin{defn}
A nearring $N$ is \textit{equiprime} if, whenever $a,x,y\in N$ and
$anx=any$ for all $n\in N$, then $a=0$ or $x=y$. A proper ideal
$I$ of $N$ is \emph{equiprime} if and only if the nearring $N/I$
is equiprime. The \emph{equiprime radical} ${\cal P}_{*}(N)$ of $N$
is defined to be the intersection of all equiprime ideals of $N$.
If no proper ideal of $N$ is equiprime, then ${\cal P}_{*}(N)$ is
defined to be $N$ itself.
\end{defn}

Note that every equiprime nearring is zero-symmetric. As stated in
\cite[Corollary~3.5]{BoothGV}, ${\cal P}_{*}$ is an ideal-hereditary
Kurosh-Amitsur radical class in the variety of zero-symmetric nearrings,
and for any zero-symmetric nearring $N$, ${\cal P}_{*}(N)$ contains
the prime radical ${\cal P}(N)$ (the intersection of the $0$-prime
ideals of $N$) \cite[Proposition~4.2]{BoothGV}.

It is known that ${\cal P}_{*}$ does not contain the (Brown McCoy)
upper radical class $\mathcal{B}_{0}$ determined by the zero-symmetric
simple nearrings with identity (this is the case even in the variety
of rings, cf.~Divinsky \cite{Div}). The problem of whether $\mathcal{B}_{0}$
contains ${\cal P}_{*}$ remained open until this point.

If $N$ is a non-equiprime zero-symmetric simple nearring with identity,
then $N$ would serve as a counter-example to ${\cal P}_{*}\subseteq{\cal B}_{0}$.
This is because the only nontrivial homomorphic image of $N$ is $N$
itself, and so it belongs to ${\cal P}_{*}$, but not to $\mathcal{B}_{0}$.

In the next two sections, we will demonstrate that a zero-symmetric
simple nearring~$N$ with identity may not be equiprime.

Before delving into that, let us make some remarks. Consider a nearring
$N$ with the additive group $(N,+)$ being a simple group; in this
case, $N$ is certainly a simple nearring. Recent constructions, as
described in \cite{KeMPW}, utilize Higman-Neumann-Neumann (HNN) extensions
to obtain infinite simple (additive) groups $G^{*}$ from a base group~$G$.
These constructions allow for the creation of zero-symmetric simple
nearrings $(N,+,\cdot)$ with identity, where $(N,+)=(G^{*},+)$.
For a comprehensive understanding of the necessary and sufficient
conditions to build a nearring with identity on a group, refer to
\cite{Clay}.

Two constructions are presented in \cite{KeMPW}. The first construction
is based on $G=(\mathbb{Z},+)$, and the second one on $G=(F(B),+)$,
where $F(B)$ represents the free group on an infinite set $B$. Surprisingly,
the second construction yields non-equiprime nearrings, while the
first construction results in an equiprime nearring! We will explore
these facts in the next section after briefly reviewing how the constructions
were carried out.

In the final section, we introduce a modification of the first construction
of \cite{KeMPW}, still based on $G=(\mathbb{Z},+)$, to obtain another
zero-symmetric simple nearring with identity that is not equiprime.

\section{Zero-symmetric simple nearrings with identity whose additive groups
are simple}

Let $G^{*}$ be an additive group. The following conditions are necessary
and sufficient for one to define a multiplication $\cdot$ on $G^{*}$
and obtain a zero-symmetric nearring with identity $(N,+,\cdot)$
where $(N,+)=(G^{*},+)$: (1) there must be a group isomorphism $f_{\gamma}:G^{*}\to H_{\gamma}$,
one for each nonzero element $\gamma\in G^{*}$, from $G^{*}$ to
a proper subgroup $H_{\gamma}$ of $G^{*}$ containing $\gamma$;
(2) for $\gamma,\zeta\in G^{*}$, $f_{\gamma}\circ f_{\zeta}=f_{f_{\gamma}(\zeta)}$;
and (3) there is an element $e\in G^{*}$ such that $f_{e}(\gamma)=f_{\gamma}(e)=\gamma$
for all $\gamma\in G^{*}$. (See Theorems~1.1 and~1.2 of \cite{Clay}.)
The nearring multiplication is then given by $\gamma\cdot\zeta=f_{\zeta}(\gamma)$
for all $\gamma,\zeta\in G^{*}$.

Let us briefly describe the construction in \cite{KeMPW} of the infinite
simple groups $G^{*}$ such that zero-symmetric simple nearrings with
identity can be built on them.

The base group $G$ is either $\mathbb{Z}=\langle1\rangle$ or a free
group of infinite rank on the basis $B$ with $1\in B$. Set $G_{-1}=\{0\}$,
$G_{0}=G$, and suppose that the group $G_{i}$ where $i\geq0$ is
already obtained.

For each distinct pair of nonzero elements $\alpha,\beta\in G_{i}$
such that either $\alpha\not\in G_{i-1}$ or $\beta\not\in G_{i-1}$,
a free generator $t_{\alpha,\beta}$ is added to $G_{i}$ with the
relation $-t_{\alpha,\beta}+\alpha+t_{\alpha,\beta}=\beta$. Denote
by ${\cal T}(G_{i})$ the set of all such $t_{\alpha,\beta}$, and
put $G_{i+1}=\langle G_{i}\cup{\cal T}(G_{i})\rangle$. Then $G_{i+1}$
is a torsion-free group where all nonzero elements in $G_{i}$ are
conjugate in $G_{i+1}$. In this way, a sequence of groups $G=G_{0}\subseteq G_{1}\subseteq G_{2}\subseteq\dots$
is obtained, and $G^{*}=\cup_{i=0}^{\infty}G_{i}$ is a simple group
since every two distinct nonzero elements $\alpha,\beta\in G^{*}$
are conjugate.

Now, for each nonzero element $\gamma\in G^{*}$, a sequence of subgroups
$H_{\gamma}=H_{\gamma,0}\subseteq H_{\gamma,1}\subseteq H_{\gamma,2}\subseteq\dots$
with $\gamma\in H_{\gamma}$ is derived such that there are isomorphisms
$f_{\gamma}^{(i)}:G_{i}\to H_{\gamma,i}$, $i=0,1,2,\dots$, having
the following properties:
\begin{enumerate}
\item $f_{\gamma}^{(0)}(1)=\gamma$,
\item for $i\geq0$, $f_{\gamma}^{(i+1)}|_{G_{i}}=f_{\gamma}^{(i)}$, and
\item $f_{\gamma}^{(i+1)}(t_{\alpha,\beta})=t_{f_{\gamma}^{(i)}(\alpha),f_{\gamma}^{(i)}(\beta)}$
for all distinct $\alpha,\beta\in G_{i}\setminus\{0\}$ with either
$\alpha\not\in G_{i-1}$ or $\beta\not\in G_{i-1}$.
\end{enumerate}
Then the subgroup $H_{\gamma}^{*}=\cup_{i=0}^{\infty}H_{\gamma,i}$
of $G^{*}$ is isomorphic to $G^{*}$ with an isomorphism $f_{\gamma}:G^{*}\to H_{\gamma}^{*}$
such that $f_{\gamma}|_{G_{i}}=f_{\gamma}^{(i)}$ for $i=0,1,2,\dots$.
Thus, if $\alpha,\beta\in G^{*}$ are distinct, then $f_{\gamma}(t_{\alpha,\beta})=t_{f_{\gamma}(\alpha),f_{\gamma}(\beta)}$.

For $\zeta,\tau\in G^{*}$, it was shown that the equality $f_{\zeta}\circ f_{\tau}=f_{f_{\zeta}(\tau)}$
holds. Therefore, one can define a binary operation $\cdot_{f}$ on
$G^{*}$, where $\alpha\cdot_{f}\beta=f_{\beta}(\alpha)$ for all
$\alpha,\beta\in G^{*}$, and obtain a zero-symmetric nearring $(N,+,\cdot)$,
where $\cdot=\cdot_{f}$, with $(N,+)=(G^{*},+)$, and such that $1$
is the identity of $N$.

Next, we will go into more detail on how the subgroups $H_{\gamma}$,
$\gamma\in G^{*}$, are constructed, and show that (1) the zero-symmetric
nearrings with identity obtained on $(F(B),+)$, with $B$ an infinite
set, are not equiprime, and (2) the zero-symmetric nearring with identity
obtained on $(\mathbb{Z},+)$ is equiprime.

\subsection{Non-equiprime zero-symmetric simple nearrings with identity}

Consider first the case when the starting group $G$ is a free group
of infinite rank on the basis~$B$ with $1\in B$.

After $G^{*}$ is constructed as described above, the basis $B$ of
the starting group $G=G_{0}$ is indexed by the smallest ordinal number
$X$ with $|X|=|B|$ such that $1$ is the first element in $B$.
So $B=\{\pi_{\iota}\mid\iota\in X\}$ and $\pi_{0}=1$.

It was shown that every nonzero $\gamma\in G^{*}$ uniquely determines
an element $\mu_{\gamma}\in X$ such that the set $B_{\gamma}=\{\gamma\}\cup\{\pi_{\mu_{\gamma}+\iota}\mid\iota\succ0\}$
is independent. Thus, the subgroup $H_{\gamma}$ generated by $B_{\gamma}$
is isomorphic to $G_{0}$ with the isomorphism $f_{\gamma}^{(0)}:G_{0}\to H_{\gamma}$
satisfying $f_{\gamma}^{(0)}(1)=\gamma$ and $f_{\gamma}^{(0)}(\pi_{\iota})=\pi_{\mu_{\gamma}+\iota}$
for $\iota\succ0$. Then $f_{\gamma}:G^{*}\to H_{\gamma}^{*}$ satisfies
$f_{\gamma}|_{G_{i}}=f_{\gamma}^{(i)}$ for $i=0,1,2,\dots$, and
$f_{\gamma}(t_{\alpha,\beta})=t_{f_{\gamma}(\alpha),f_{\gamma}(\beta)}$
for all distinct $\alpha,\beta\in G^{*}\setminus\{0\}$. The nearring
$(N,+,\cdot)$ with $a\cdot b=f_{b}(a)$ for all $a,b\in N$ is a
zero-symmetric simple nearring with identity $1$ and $(N,+)=(G^{*},+)$.

Of particular relevance for us here is that $\mu_{k\pi_{\iota}}=\iota$
for all $\iota\in X$ and $k\in\mathbb{Z}\setminus\{0\}$, and $\mu_{f_{\gamma}(x)}=\mu_{\gamma}+\mu_{x}$
for all $x\in G^{*}\setminus\{0\}$.

We claim that this nearring is not equiprime. To see this, take $a=b=\pi_{1}$
and $c=2\pi_{1}$. Then $a,b,c\in G_{0}$ and $\mu_{a}=\mu_{b}=\mu_{c}=1$.
Let $x\in G^{*}\setminus\{0\}$. From
\[
\mu_{f_{b}(x)}=\pi_{\mu_{b}+\mu_{x}}=\pi_{\mu_{c}+\mu_{x}}=\mu_{f_{c}(x)},
\]
we get
\begin{align*}
axb & =f_{f_{b}(x)}(a)=f_{f_{b}(x)}(\pi_{1})=f_{f_{b}(x)}^{(0)}(\pi_{1})=\pi_{\mu_{f_{b}(x)}+1}\\
 & =\pi_{\mu_{f_{c}(x)}+1}=f_{f_{c}(x)}^{(0)}(\pi_{1})=f_{f_{c}(x)}(\pi_{1})=f_{f_{c}(x)}(a)=axc.
\end{align*}
Clearly, $\pi_{1}0b=0=\pi_{1}0c$. Since $b\not=c$, $N$ is not equiprime.

Now, being a nonequiprime simple nearring with identity, $N$ has
to have a nontrivial invariant subgroup as indicated in~\cite{Veldsman}.
This means that there is a nontrivial proper additive subgroup $M$
of $N$ such that and $b\cdot M\subseteq M$ for all $b\in N$ ($M$
is left invariant), and $M\cdot b\subseteq M$ for all $b\in N$ ($M$
is right invariant).

Let $W_{0}=\gen[\{\pi_{\iota}\mid\iota\succ0\}]$. For each $\ell\geq1$,
we collect those elements $t_{\alpha,\beta}\in{\cal T}(G_{\ell-1})$
with $\alpha,\beta\in W_{\ell-1}$, and put them in $T_{\ell-1}$,
then set $W_{\ell}=\gen[W_{\ell-1}\cup T_{\ell-1})]\leq G_{\ell+1}$.
The union of the sequence $W_{0}\subseteq W_{1}\subseteq W_{2}\subseteq\dots$
of subgroups of $G^{*}$ results in a nontrivial proper subgroup $W^{*}=\cup_{\ell=0}^{\infty}W_{\ell}$
of $G^{*}$. Set $W_{-1}=\{0\}$.

We claim that $W^{*}$ is an invariant subgroup of $N$. Let $\zeta\in N$.

We first show that $\zeta\cdot W^{*}\subseteq W^{*}$. Take any $\tau\in W^{*}\setminus\{0\}$.
Then $\tau\in W_{i}\setminus W_{i-1}$ for some $i\geq0$. We have
$H_{\tau,0}=\gen[\{\tau\}\cup\{\pi_{\iota}\mid\iota\succ\mu_{\tau}\}]\subseteq W_{0}$,
$H_{\tau,1}=\gen[H_{\tau,0}\cup T_{\tau,0}]\subseteq W_{1}$, $\dots$,
$H_{\tau,n+1}=\gen[H_{\tau,n}\cup T_{\tau,n}]\subseteq W_{n+1}$,
etc. Hence $H_{\tau}^{*}=\cup_{n=0}^{\infty}H_{\tau,n}\subseteq\cup_{n=0}^{\infty}W_{n}=W^{*}$.
From $f_{\tau}(G^{*})=H_{\tau}^{*}$, we see that $\zeta\cdot\tau=f_{\tau}(\zeta)\in H_{\tau}^{*}\subseteq W^{*}$.
Therefore, $\zeta\cdot W^{*}\subseteq W^{*}$. Consequently, $W^{*}$
is left invariant.

Next, we will show that $f_{\zeta}(W^{*})\subseteq W^{*}$. When this
is done, it would follow that $\tau\cdot\zeta=f_{\zeta}(\tau)\in W^{*}$
for all $\tau\in W^{*}$ and so $W^{*}\cdot N\subseteq W^{*}$. For
any $\iota\in X$ and $\iota\succ0$, we have
\begin{equation}
f_{\zeta}(\pi_{\iota})=f_{\zeta}^{(0)}(\pi_{\iota})=\pi_{\mu_{\zeta}+\iota}\in W_{0},\label{eq:base-1}
\end{equation}
and so $f_{\zeta}(W_{0})\subseteq W_{0}$. For $t_{\alpha,\beta}\in{\cal T}(G_{0})$
with $\alpha,\beta\in W_{0}\subseteq G_{0}$, we have
\[
f_{\zeta}(t_{\alpha,\beta})=t_{f_{\zeta}^{(0)}(\alpha),f_{\zeta}^{(0)}(\beta)}\in W_{1}.
\]
Therefore, $f_{\zeta}(W_{1})\subseteq W_{1}$. In the same manner,
we see that $f_{\zeta}(W_{2})\subseteq W_{2}$, $f_{\zeta}(W_{3})\subseteq W_{3}$,
and so on. Therefore, $f_{\zeta}(W^{*})\subseteq W^{*}$. We are done.

\subsection{An equiprime zero-symmetric simple nearring with identity}

Next, let $(N,+,\cdot)$ be the nearring with $(N,+)=(G^{*},+)$ where
$G=\mathbb{Z}$. Let us show that this nearring is equiprime. We need
the following lemma.
\begin{lem}
\label{lem:Z}Let $a,b,c\in G^{*}$ with $a\not=0$. If $f_{t_{b,-b}}(a)=f_{t_{c,-c}}(a)$,
then $b=c$.
\end{lem}

\begin{proof}
Let $G_{i}'=G_{i}\setminus G_{i-1}$, $i=0,1,2,\dots,$ then $G^{*}=\{0\}\cup G_{0}'\cup G_{1}'\cup\dots$,
where the union is disjoint, and $a\in G_{i}'$ for some unique $i\geq0$.
We will use induction on $i$ to show that $b=c$. Suppose that $a\in G_{0}'=\mathbb{Z}\setminus\{0\}$.
Then
\[
f_{t_{b,-b}}(a)=\epsilon\big(\underbrace{f_{t_{b,-b}}(1)+\dots+f_{t_{b,-b}}(1)}_{|a|\text{ times}}\big)=\epsilon\big(\underbrace{t_{b,-b}+\dots+t_{b,-b}}_{|a|\text{ times}}\big)
\]
where $\epsilon=\frac{a}{|a|}$. Similarly,
\[
f_{t_{c,-c}}(a)=\epsilon\big(\underbrace{t_{c,-c}+\dots+t_{c,-c}}_{|a|\text{ times}}\big).
\]
If $b\not=c$, then $t_{c,-c}$ and $t_{b,-b}$ are distinct free
generators having no relation. But then
\[
f_{t_{b,-b}}(a)=\epsilon(\underbrace{t_{b,-b}+\dots+t_{b,-b}}_{|a|\text{ times}})\not=\epsilon(\underbrace{t_{c,-c}+\dots+t_{c,-c}}_{|a|\text{ times}})=f_{t_{c,-c}}(a),
\]
a contradiction. Therefore, we have $b=c$ in this case.

Assume that $i>0$, $a\in G_{i}'$, and if there is some $e\in G_{j}'$,
where $0\leq j<i$, such that $f_{t_{b,-b}}(e)=f_{t_{c,-c}}(e)$,
then $b=c$. The HNN construction of $G^{*}$ gives $G_{i}=G_{i-1}*{\cal T}_{i}$,
the free product of $G_{i-1}$ and the free group ${\cal T}_{i}$
generated by $t_{\alpha,\beta}$ where $\alpha,\beta\in G_{i-1}\setminus\{0\}$
are distinct with either $\alpha\not\in G_{i-2}$ or $\beta\not\in G_{i-2}$.
Thus, $a$ can be written uniquely as $a=\gamma_{0}+t_{1}+\gamma_{1}+\dots+t_{k}+\gamma_{k}$,
where $k\geq1$ with $\gamma_{0},\gamma_{1},\dots,\gamma_{k}\in G_{i-1}$
and each $t_{1},t_{2},\dots,t_{k}\in{\cal T}_{i}$. Then
\[
f_{t_{b,-b}}(a)=f_{t_{b,-b}}(\gamma_{0})+f_{t_{b,-b}}(t_{1})+f_{t_{b,-b}}(\gamma_{1})+\dots+f_{t_{b,-b}}(t_{k})+f_{t_{b,-b}}(\gamma_{k})
\]
and
\[
f_{t_{c,-c}}(a)=f_{t_{c,-c}}(\gamma_{0})+f_{t_{c,-c}}(t_{1})+f_{t_{t_{c,-c}}}(\gamma_{1})+\dots+f_{t_{t_{c,-c}}}(t_{k})+f_{t_{t_{c,-c}}}(\gamma_{k}).
\]
are the unique expressions of $f_{t_{b,-b}}(a)$ and $f_{t_{c,-c}}(a)$
in some $G_{s}=G_{s-1}*{\cal T}_{s}$, where $s\geq i$. In particular,
$f_{t_{b,-b}}(t_{1})=f_{t_{c,-c}}(t_{1})$. As $t_{1}\in{\cal T}_{i}$,
we have $t_{1}=\sum_{j=1}^{\ell}n_{j}t_{\alpha_{j},\beta_{j}}$ for
some $n_{j}\in\mathbb{Z}\setminus\{0\}$ and $\alpha_{1},\dots,\alpha_{\ell},\beta_{1},\dots,\beta_{\ell}\in G_{i-1}$
with $\alpha_{j}\not=\beta_{j}$ for all $j$. From $f_{t_{b,-b}}(a)=f_{t_{c,-c}}(a)$
and that $f_{t_{b,-b}}(t_{\alpha_{1},\beta_{1}})=t_{f_{t_{b,-b}}(\alpha_{1}),f_{t_{b,-b}}(\beta_{1})}$
and $f_{t_{c,-c}}(t_{\alpha_{1},\beta_{1}})=t_{f_{t_{c,-c}}(\alpha_{1}),f_{t_{c,-c}}(\beta_{1})}$
are free generators in ${\cal T}_{s}$, we deduce that $t_{f_{t_{b,-b}}(\alpha_{1}),f_{t_{b,-b}}(\beta_{1})}=t_{f_{t_{c,-c}}(\alpha_{1}),f_{t_{c,-c}}(\beta_{1})}$.
Consequently, $f_{t_{b,-b}}(\alpha_{1})=f_{t_{c,-c}}(\alpha_{1})$.
It follows from the induction hypothesis that $b=c$, and we are done.
\end{proof}
Assume that $a,b,c\in N$ with $a\not=0$, and $axb=axc$ for all
$x\in N$. In particular, $a\cdot t_{1,-1}\cdot b=a\cdot t_{1,-1}\cdot c$.
We have $t_{1,-1}\cdot b=f_{b}(t_{1,-1})=t_{f_{b}(1),f_{b}(-1)}=t_{b,-b}$.
Likewise, $t_{1,-1}\cdot c=t_{c,-c}$. Thus, $f_{t_{b,-b}}(a)=a\cdot t_{b,-b}=a\cdot t_{c,-c}=f_{t_{c,-c}}(a)$,
and so $b=c$ by Lemma~\ref{lem:Z}. Therefore, $N$ is an equiprime
nearring.

\section{Another nonequiprime example}

Here we present an alternative construction using HNN extensions,
beginning with $G=(\mathbb{Z},+)$, to produce a new simple group
$G^{*}$. This group $G^{*}$ serves as the additive group of a countable
zero-symmetric nearring with identity, where the nearring is not equiprime.

Put $G_{0}=G$. Let $G_{-1}=\{0\}$. We want to get a sequence $G=G_{0}\subseteq G_{1}\subseteq G_{2}\subseteq\dots$
such that every pair of distinct nonzero elements in $G_{i}$ is conjugate
in $G_{i+1}$ for $i=0,1,2,\dots$. Suppose that we have already obtained
$G_{i}$. Then $G_{i+1}$ is given by the following construction.

Let
\[
{\cal T}(G_{i})=\{t_{\alpha,\beta}\mid\alpha,\beta\in G_{i}\setminus\{0\},\text{ }\alpha\not=\beta,\text{ }\alpha\not\in G_{i-1}\text{ or }\beta\not\in G_{i-1}\},
\]
where each $t_{\alpha,\beta}$ is a free generator, subject only to
the relation $-t_{\alpha,\beta}+\alpha+t_{\alpha,\beta}=\beta$. Also,
take a free generator $\varpi_{i}\not\in\left\langle G_{i}\cup{\cal T}(G_{i})\right\rangle $
such that $\varpi_{i}$ commutes with every element in $G_{i}\cup{\cal T}(G_{i})$.
Then we put put $G_{i+1}=\gen[G_{i}\cup{\cal T}(G_{i})\cup\{\varpi_{i}\}]=\gen[G_{i}\cup{\cal T}(G_{i})]\oplus\gen[\varpi_{i}]$.

Now, the group $G^{*}=\cup_{i=0}^{\infty}G_{i}$ is a countable torsion-free
group, and every pair of distinct nonzero elements in it are conjugate.

Let $\zeta\in G^{*}\setminus\{0\}$. Then $\zeta\in G_{i}\setminus G_{i-1}$
for some unique $i\geq0$, which we will denote by $|\zeta|=i$. We
want to associate with $\zeta$ a subgroup $H_{\zeta}^{*}$ such that
there is an isomorphism $f_{\zeta}:G^{*}\to H_{\zeta}$ with $f_{\zeta}(1)=\zeta$.
This will be done by building up a sequence $H_{\zeta,0}\subseteq H_{\zeta,1}\subseteq H_{\zeta,2}\subseteq\dots$
of subgroups of $G^{*}$ such that each $H_{\zeta,j}\leq G_{i+j}$
with an isomorphisms $f_{\zeta}^{(j)}:G_{j}\to H_{\zeta,j}$, for
$j=0,1,2,\dots$. Then the subgroup $H_{\zeta}^{*}=\cup_{j=0}^{\infty}H_{\zeta,j}$
of $G^{*}$ is isomorphic to $G^{*}$ with an isomorphism $f_{\zeta}:G^{*}\to H_{\zeta}^{*}$
satisfying $f_{\zeta}|_{G_{j}}=f_{\zeta}^{(j)}$ for $j=0,1,2,\dots$.
We shall proceed by induction on $j$, and, for the start, set $H_{\zeta,0}=\gen[\zeta]\leq G_{i}$.
The map $f_{\zeta}^{(0)}:G_{0}\to H_{\zeta,0}$ given by $f_{\zeta}^{(0)}(m)=m\zeta$
for all $m\in\mathbb{Z}$ is a group isomorphism.

Assume that we have already obtained $H_{\zeta,j}\leq G_{i+j}$ and
$f_{\zeta}^{(j)}:G_{j}\to H_{\zeta,j}$. We collect those $t_{f_{\zeta}^{(j)}(\alpha),f_{\zeta}^{(j)}(\beta)}\in T_{0}(G_{i+j})\subseteq G_{i+(j+1)}$,
where $\alpha,\beta\in G_{j}\setminus\{0\}$ are distinct, and put
them into the set $\widetilde{{\cal T}_{j}}(\zeta)$. Then define
\[
H_{\zeta,j+1}=\gen[H_{\zeta,j}\cup\widetilde{{\cal T}_{j}}(\zeta)\cup\{\varpi_{i+j)}\}]=\gen[H_{\zeta,j}\cup\widetilde{{\cal T}_{j}}(\zeta)]\oplus\gen[\varpi_{i+j}]
\]
which is a subgroup of $G_{i+(j+1)}$. Then $H_{\zeta,j+1}$ is isomorphic
to $G_{j+1}=\gen[G_{j}\cup{\cal T}(G_{j})]\oplus\gen[\varpi_{j}]$
with an isomorphism $f_{\zeta}^{(j+1)}:G_{j+1}\to H_{\zeta,j+1}$
determined by
\begin{enumerate}
\item $f_{\zeta}^{(j+1)}(\alpha)=f_{\zeta}^{(j)}(\alpha)$ if $\alpha\in G_{j}$,
\item $f_{\zeta}^{(j+1)}(t_{\alpha,\beta})=t_{f_{\zeta}^{(j)}(\alpha),f_{\zeta}^{(j)}(\beta)}$
if $t_{\alpha,\beta}\in T_{0}(G_{j})$, where $\alpha,\beta\in G_{j}\setminus\{0\}$
are distinct,
\item $f_{\zeta}^{(j+1)}(\varpi_{j})=\varpi_{i+j}$.
\end{enumerate}
This completes the inductive step, and yields the desire subgroup
$H_{\zeta}^{*}=\cup_{j=0}^{\infty}H_{\zeta,j}$ of $G^{*}$ with the
isomorphism $f_{\zeta}:G^{*}\to H_{\zeta}^{*}$. Here, we have that
the value $f_{\zeta}(1)=\zeta$ determines $f_{\zeta}^{(0)}$, $f_{\zeta}^{(0)}$
determines $f_{\zeta}^{(1)}$, $f_{\zeta}^{(1)}$ determines $f_{\zeta}^{(2)}$,
etc. Moreover, we have, for each $j\geq0$,
\begin{equation}
f_{\zeta}(\varpi_{j})=f_{\zeta}^{(j+1)}(\varpi_{j})=\varpi_{i+j}=\varpi_{|\zeta|+j}.\label{eq:pi_j}
\end{equation}

Now, for $\zeta,\tau\in G^{*}$, from
\[
(f_{\zeta}\circ f_{\tau})(1)=f_{\zeta}(f_{\tau}(1))=f_{\zeta}(\tau)=f_{f_{\zeta}(\tau)}(1),
\]
we infer that $f_{\zeta}\circ f_{\tau}=f_{f_{\zeta}(\tau)}$. Therefore,
by~\cite[Theorems 1.1 and 1.2]{Clay}, with the multiplication $\cdot$
on $G^{*}$ given by $\zeta\cdot\tau=f_{\tau}(\zeta)$ for all $\zeta,\tau\in G^{*}$,
we have a zero-symmetric simple nearring $N=(N,+,\cdot)$ with identity,
having $(N,+)=(G^{*},+)$. Here, the nearring identity is~$1$.

We claim that $N$ is not equiprime. To see this, we first make an
observation.

Let $\zeta\in G_{0}\setminus\{0,1\}=\mathbb{Z}\setminus\{0,1\}$ and
$\lambda\in G^{*}$. With $|\zeta|=0$, by (\ref{eq:pi_j}), we have
$f_{\zeta}(\varpi_{j})=\varpi_{|\zeta|+j}=\varpi_{j}$ for $j\geq0$.
Now, if $\lambda\in G_{0}=\mathbb{Z}$, then $f_{\zeta}(\lambda)=\lambda f_{\zeta}(1)=\lambda\zeta\in G_{0}$.
Suppose that $k\geq0$, $\lambda\in G_{k+1}$, and that $f_{\zeta}(\tau)\in G_{k+1}$
for all $\tau\in G_{j}$ with $j\leq k$. As $G_{k+1}=\gen[G_{k}\cup{\cal T}(G_{k})\cup\{\varpi_{k}\}]$,
we can write
\[
\lambda=m\varpi_{k}+\sum_{\ell=1}^{s}\big(\epsilon_{\ell}\lambda_{\ell}+\epsilon_{\ell}'t_{\alpha_{\ell},\beta_{\ell}}\big)
\]
for some $m,\epsilon_{\ell},\epsilon_{\ell}'\in\mathbb{Z}$ and $\lambda_{\ell},\alpha_{\ell},\beta_{\ell}\in G_{k}$
($\ell=1,2,\dots,s$). Hence
\[
f_{\zeta}(\lambda)=m\varpi_{k}+\sum_{\ell=1}^{s}\big(\epsilon_{\ell}f_{\zeta}(\lambda_{\ell})+\epsilon_{\ell}'t_{f_{\zeta}(\alpha_{\ell}),f_{\zeta}(\beta_{\ell})}\big)
\]
which is an element in $G_{k+1}$. Thus, by induction, we have $f_{\zeta}(\lambda)\in G_{|\lambda|}$.
Consequently,
\begin{equation}
|f_{\zeta}(\lambda)|=|\lambda|.\label{eq:same_layer}
\end{equation}

Now, take $\zeta_{1},\zeta_{2}\in G_{0}\setminus\{0,1\}$ with $\zeta_{1}\not=\zeta_{2}$.
We have
\[
\varpi_{0}\cdot0\cdot\zeta_{1}=0=\varpi_{0}\cdot0\cdot\zeta_{2},
\]
and for every $\tau\in G^{*}\setminus\{0\}$,
\[
\varpi_{0}\cdot\tau\cdot\zeta_{1}=f_{f_{\zeta_{1}}(\tau)}(\varpi_{0})=\varpi_{|f_{\zeta_{1}}(\tau)|}=\varpi_{|\tau|}=\varpi_{|f_{\zeta_{2}}(\tau)|}=f_{f_{\zeta_{2}}(\tau)}(\varpi_{0})=\varpi_{0}\cdot\tau\cdot\zeta_{2}
\]
by (\ref{eq:pi_j}) and (\ref{eq:same_layer}). As $\varpi_{0}\not=0$
and $\zeta_{1}\not=\zeta_{2}$, this shows that $N$ is not equiprime.
\smallskip{}

We can also find a nontrivial invariant subgroup of $N$.

Note that for any $\zeta\in G^{*}\setminus\{0,1\}$, the proper subgroup
$H_{\zeta}^{*}$ of $G^{*}$ is left invariant in $N$. To see this,
we let $\tau\in H_{\zeta}^{*}\setminus\{0\}$ and put $H_{\zeta,-1}=\{0\}$;
hence $\tau\in H_{\zeta,i}\setminus H_{\zeta,i-1}$ for some $i\geq0$.
Let $K_{0}=H_{\tau,0}=\gen[\tau]$. Certainly, $K_{0}\subseteq H_{\zeta,i}$.
Then, $K_{1}=\gen[K_{0}\cup\widetilde{{\cal T}_{0}}(\tau)\cup\{\varpi_{i}\}]\subseteq H_{\zeta,i+1}$,
$\dots$, $K_{j+1}=\gen[K_{j}\cup\widetilde{{\cal T}_{j}}(\tau)\cup\{\varpi_{i+j}\}]\subseteq H_{\zeta,i+j+1}$,
etc. Hence $K^{*}=\cup_{j=0}^{\infty}K_{j}\subseteq\cup_{j=0}^{\infty}H_{\zeta,i+j}\subseteq H_{\zeta}^{*}$.
From $f_{\tau}:G^{*}\to K^{*}\subseteq H^{*}$, we see that $\gamma\cdot\tau=f_{\tau}(\gamma)\in H_{\zeta}^{*}$
for all $\gamma\in N$. Thus, $\gamma\cdot H_{\zeta}^{*}\subseteq H_{\zeta}^{*}$
for all $\gamma\in N$, and so $H_{\zeta}^{*}$ is left invariant.
In particular, if we take $\zeta=\varpi_{0}$, then the proper subgroup
$H_{\varpi_{0}}^{*}$ of $G^{*}$ is a left invariant subgroup of~$N$.

We claim that $H_{\varpi_{0}}^{*}$ is also right invariant. Take
an arbitrary $\gamma\in G^{*}\setminus\{0\}$. We have to show that
$H_{\varpi_{0}}^{*}\cdot\gamma\subseteq H_{\varpi_{0}}^{*}$, i.e.,
$f_{\gamma}(\tau)=\tau\cdot\gamma\in H_{\varpi_{0}}^{*}$ for all
$\tau\in H_{\varpi_{0}}^{*}$.

Recall that $H_{\varpi_{0}}^{*}=\cup_{j=0}^{\infty}H_{\varpi_{0},j}$
with $H_{\varpi_{0},0}=\gen[\varpi_{0}]$, and
\[
H_{\varpi_{0},j+1}=\gen[H_{\varpi_{0},j}\cup\widetilde{{\cal T}_{j}}(\varpi_{0})\cup\{\varpi_{j}\}],\text{ }j\geq0.
\]
Here for each $j\geq0$,
\[
\widetilde{{\cal T}_{j}}(\varpi_{0})=\{t_{f_{\varpi_{0}}^{(j)}(\alpha),f_{\varpi_{0}}^{(j)}(\beta)}\in H_{\varpi_{0},j+1}\mid\alpha,\beta\in G_{0}\setminus\{0\}\text{ with }\alpha\not=\beta\}.
\]

Assume that $\gamma\in G_{i}\setminus G_{i-1}$. From (\ref{eq:pi_j}),
we have
\begin{equation}
f_{\gamma}(\varpi_{0})=f_{\gamma}^{(1)}(\varpi_{0})=\varpi_{i},\label{eq:base}
\end{equation}
and so $f_{\gamma}(H_{\varpi_{0},0})=f_{\gamma}^{(1)}(H_{\varpi_{0},0})\subseteq H_{\varpi_{0},i}$.
For $\alpha,\beta\in H_{\varpi_{0},0}\setminus\{0\}$, $\alpha\not=\beta$,
we have $t_{\alpha,\beta}\in H_{\varpi_{0},1}$ and so
\[
f_{\gamma}(t_{\alpha,\beta})=t_{f_{\gamma}^{(1)}(\alpha),f_{\gamma}^{(1)}(\beta)}\in H_{\varpi_{0},i+1}.
\]
Since $f_{\gamma}(\varpi_{1})=\varpi_{i+1}$, we see that $f_{\gamma}(H_{\varpi_{0},1})\subseteq H_{\varpi_{0},i+1}$.
In the same manner, we see that $f_{\gamma}(H_{\varpi_{0},2})\subseteq H_{\varpi_{0},i+2}$,
$f_{\gamma}(H_{\varpi_{0},3})\subseteq H_{\varpi_{0},i+3}$, and so
on. Therefore, $f_{\gamma}(H_{\varpi_{0}}^{*})\subseteq H_{\varpi_{0}}^{*}$.
We are done. \smallskip{}

We close our discussions with the following remark.
\begin{rem}
When $N$ is a zero-symmetric simple nearring with identity, the $k\times k$
matrix nearring $M_{k}(N)$ (for any $k\ge1$) is also a simple zero-symmetric
nearring with identity. See, for example, \cite{Mey}. It was shown
in \cite{Veldsman} that $N$ is equiprime if and only if $M_{k}(N)$
is equiprime. Consequently, for any of the foregoing nearrings $N$
that turned out to be simple, zero-symmetric with identity, and not
equiprime, the matrix nearring $M_{k}(N)$ is also an example of a
simple, zero-symmetric nearring with identity that is not equiprime.
\end{rem}

\section*{Acknowledgment}

The second author would like to thank prof. Ke for his assistance
and hospitality during his stay in Tainan. He would also like to thank
the University of the Free State and the Research Center for Theoretical
Sciences in the College of Sciences, National Cheng Kung University
for financial assistance.

\vspace*{1pc}
\noindent\textbf{Funding.} No funding was received to assist with the preparation of this manuscript.\\
\noindent\textbf{Competing Interests.} The authors have no relevant financial or non-financial interests to disclose.\\
\noindent\textbf{Conflict of interest.} The authors declare that
there are no conflicts of interests.\\

\end{document}